\documentclass[twosided,reqno]{amsart}

\usepackage{amsmath}
\usepackage{amssymb}
\usepackage{amsthm}

\newcommand{\iz}{{\int_{{\mathbb{Z}}_p}}}

\theoremstyle{theorem}
\newtheorem{theorem}{\scshape Theorem }[section]

\theoremstyle{definition}

\numberwithin{equation}{section}

\begin{document}

\title[On degenerate  $q$-Euler polynomials]{On degenerate  $q$-Euler polynomials}

\author{Dmitry V. Dolgy$^1$}
\address{$^1$ Institute of Natural Sciences, Far Eastern Federal University,
Vladivostok, 690950, Russia.}

\author{Taekyun Kim$^2$}
\address{$^2$ Department of Mathematics, Kwangwoon University, Seoul 139-701, Republic of Korea.}
\email{tkkim@kw.ac.kr}

\author{Jin-Woo Park$^{3,*}$}
\address{$^3$ Department of Mathematics, Kwangwoon University, Seoul 139-701, Republic of Korea.}
\email{a0417001@knu.ac.kr}
\thanks{$*$ corresponding author}

\author{Jong-Jin Seo$^4$}
\address{$^4$ Department of Applied Mathematics, Pukyong National University,
Pusan, Republic of Korea}
\email{seo2011@pknu.ac.kr}


\maketitle

\begin{abstract}
Recently, Kim introduced Carlitz's type $q$-Euler numbers and polynomials(see \cite{05}). In this paper, we construct the degenerate Carlitz's type $q$-Euler numbers and polynomials which are derived from the fermionic $p$-adic $q$-integral on ${\mathbb{Z}}_p$. Finally, we give some identities and properties of these numbers and polynomials
\end{abstract}

\section{Introduction}

Let $p$ be a fixed odd prime number. Throughout this paper, ${\mathbb{Z}}_p$, ${\mathbb{Q}}_p$ and ${\mathbb{C}}_p$ will, respectively, denote the ring of $p$-adic rational integers, the field of $p$-adic rational numbers and the completion of the algebraic closure of ${\mathbb{Q}}_p$. Let $\nu_p$ be the normalized exponential valuation of ${\mathbb{C}}_p$ with $|p|_p=p^{-\nu_p(p)}=\frac{1}{p}$. Let $q$ be an indeterminate in ${\mathbb{C}}_p$ such that $|q-1|_p<p^{-\frac{1}{p-1}}$ so that $q^x=\exp(x\log q)$. The {\it{$q$-analogue of $x$}} is defined as $[x]_q=\frac{1-q^x}{1-q}$. Note that $\lim_{q\rightarrow 1}=x$.

As is well known, the {\it{Euler polynomials}} are defined by the generating function to be
\begin{equation}\label{1}
\frac{2}{e^t+1}e^{xt}=\sum_{n=0} ^{\infty}E_n(x)\frac{t^n}{n!},{\text{ (see [1-11])}}.
\end{equation}
When $x=0$, $E_n=E_n(0)$ are called the {\it{Euler numbers}}.  From \eqref{1}, we can derive the recurrence relation for the Euler numbers as follows:
\begin{equation}\label{2}
(E+1)^n+E_n=2\delta_{0,n},~(n\geq1),
\end{equation}
with the usual convention about replacing $E^n$ by $E_n$ (see [1-12]).

In \cite{05,06}, Kim considered Carlitz's type $q$-Euler numbers as follows:
\begin{equation}\label{3}
E_{0,q}=1,~q(qE_q+1)^n+E_{n,q}=0,~(n\geq1),
\end{equation}
with the usual convention about replacing $E_q ^n$ by $E_{n,q}$.

The $q$-Euler polynomials are defined as
\begin{equation}\label{4}
E_{n,q} (x)=\sum_{l=0} ^n \binom{n}{l}q^{lx}E_{l,q}[x]_q ^{n-l},~(n\geq0),{\text{ (see \cite{05,06})}}.
\end{equation}

Let $C({\mathbb{Z}}_p)$ be the space of all continuous functions on ${\mathbb{Z}}_p$. For $f\in C({\mathbb{Z}}_p)$, the {\it{fermionic $p$-adic $q$-integral on ${\mathbb{Z}}_p$}} is defined by Kim to be
\begin{equation}\label{5}
I_{-q}(f)=\iz f(x)d\mu_{-q}(x)=\lim_{N\rightarrow \infty}\frac{1}{\left[p^N\right]_{-q}}\sum_{x=0} ^{p^N-1}f(x)(-q)^x,
\end{equation}
where $[x]_q=\frac{1-(-q)^x}{1+q}$, (see \cite{05,06}).

From \eqref{5}, we have
\begin{equation*}
qI_{-q}(f_1)+I_{-q}(f)=[2]_qf(0),{\text{ where }}f_1(x)=f(x+1).
\end{equation*}
The $q$-Euler polynomials can be represented by the fermionic $p$-adic $q$-integral on ${\mathbb{Z}}_p$ as follows:
\begin{equation}\label{6}
\sum_{n=0} ^{\infty}E_{n,q}(x)\frac{t^n}{n!}=\iz e^{[x+y]_qt}d\mu_{-q}(y),{\text{ (see \cite{05,06})}}.
\end{equation}
Thus, by \eqref{6}, we get
\begin{equation}\label{7}
\iz[x+y]_q ^nd\mu_{-q}(y)=E_{n,q}(x),~(n\geq 0).
\end{equation}

In this paper, we consider the degenerate Carlitz $q$-Euler polynomials and numbers which are derived from the fermionic $p$-adic $q$-integral on ${\mathbb{Z}}_p$ and we investigate some properties and identities of those polynomials.

\section{Degenerate Carlitz $q$-Bernoulli numbers and polynomials}

In this section, we assume that $\lambda,t\in{\mathbb{C}}_p$ with $|\lambda t|_p<p^{-\frac{1}{p-1}}$ so that $(1+\lambda t)^{\frac{x}{\lambda}}=\exp\left(\frac{x}{\lambda}\log(1+\lambda t)\right)$.

It is well known that
\begin{equation}\label{8}
\lim_{\lambda\rightarrow 0}(1+\lambda t)^{\frac{1}{\lambda}}=\lim_{\lambda\rightarrow0}\sum_{n=0} ^{\infty}\binom{\frac{1}{\lambda}}{n}\lambda^nt^n=\sum_{n=0} ^{\infty}\frac{t^n}{n!}=e^{xt}.
\end{equation}
From \eqref{6} and \eqref{8}, we consider the {\it{degenerate $q$-Euler polynomials}} which are given by the generating function to be

\begin{equation}\label{9}
\iz(1+\lambda t)^{\frac{1}{\lambda}[x+y]_q}d\mu_{-q}(y)=\sum_{n=0} ^{\infty} E_{n,q}(x|\lambda)\frac{t^n}{n!}.
\end{equation}
Note that
\begin{equation}\label{10}
\begin{split}
\sum_{n=0} ^{\infty}\lim_{\lambda\rightarrow0}E_{n,q} (x|\lambda)\frac{t^n}{n!}=&\lim_{\lambda\rightarrow0}\iz\left(1+\lambda t\right)^{\frac{[x+y]_q}{\lambda}}d\mu_{-q}(y)\\
=&\iz e^{[x+y]_qt}d\mu_{-q  }(y)=\sum_{n=0} ^{\infty}E_{n,q}(x)\frac{t^n}{n!}.
\end{split}
\end{equation}
Thus, by \eqref{10}, we get
\begin{equation}\label{11}
\lim_{\lambda\rightarrow0}E_{n,q}(x|\lambda)=E_{n,q}(x),~(n\geq0).
\end{equation}
When $x=0$, $E_{n,q}(\lambda)=E_{n,q}(0|\lambda)$ are called {\it{degenerate $q$-Euler numbers}}.

From \eqref{9}, we have
\begin{equation}\label{12}
\begin{split}
\sum_{n=0} ^{\infty}E_{n,q}(x|\lambda)\frac{t^n}{n!}=&\sum_{m=0} ^{\infty}\lambda^{-m}\iz[x+y]_q ^md\mu_{-q}(y)\frac{\left(\log(1+t)\right)^m}{m!}\\
=&\sum_{m=0} ^{\infty}\lambda^{-m}E_{m,q}(x)\sum_{n=m} ^{\infty}S_1(n,m)\frac{\lambda^nt^n}{n!}\\
=&\sum_{n=0} ^{\infty}\left(\sum_{m=0} ^n\lambda^{n-m}E_{m,q}(x)S_1(n,m)\right)\frac{t^n}{n!}.
\end{split}
\end{equation}
By comparing the coefficients on the both sides of \eqref{12}, we obtain the following theorem.
\begin{theorem}\label{thm1}
For $n \geq 0$, we have
\begin{equation*}
E_{n,q}(x|\lambda)=\sum_{m=0} ^n \lambda^{n-m}E_{m,q}(x)S_1(n,m),.
\end{equation*}
where $S_1(n,m)$ is the Stirling number of the first kind.
\end{theorem}
Replacing $t$ by $\frac{1}{\lambda}\left(e^{\lambda t}-1\right)$ in \eqref{9}, we get
\begin{equation}\label{13}
\begin{split}
\sum_{m=0} ^{\infty}E_{m,q}(x|\lambda)\lambda^{-m}\frac{\left(e^{\lambda t}-1\right)^m}{m!}=&\iz e^{[x+y]_qt}d\mu_{-q}(y)\\
=&\sum_{n=0} ^{\infty}E_{n,q}(x)\frac{t^n}{n!}.
\end{split}
\end{equation}
On the other hand,
\begin{equation}\label{14}
\begin{split}
\sum_{m=0} ^{\infty}E_{m,q}(x|\lambda)\lambda^{-m}\frac{1}{m!}\left(e^{\lambda t}-1\right)^m=&\sum_{m=0} ^{\infty}E_{m,q}(x|\lambda)\lambda^{-m}\sum_{n=m} ^{\infty}S_2(n,m)\frac{\lambda^nt^n}{n!}\\
=&\sum_{n=0} ^{\infty}\left(\sum_{m=0} ^nE_{m,q}(x|\lambda)\lambda^{n-m}S_2(n,m)\right)\frac{t^n}{n!}.
\end{split}
\end{equation}
Therefore, by \eqref{13} and \eqref{14}, we obtain the following theorem.

\begin{theorem}\label{thm2}
For $n\geq 0$, we have
\begin{equation*}
E_{n,q}(x)=\sum_{m=0} ^n E_{m,q}(x|\lambda)\lambda^{n-m}S_2(n,m).
\end{equation*}
where $S_2(n,m)$ is the Stirling number of the second kind.
\end{theorem}

We observe that
\begin{equation}\label{15}
(1+\lambda t)^{\frac{[x+y]_q}{\lambda}}=\sum_{n=0} ^{\infty}\left(\frac{[x+y]_q}{\lambda}\right)_n\frac{\lambda^nt^n}{n!}
\end{equation}
where $(x)_n=x(x-1)\cdots(x-n+1)$.

Now we define the {\it{degenerate falling factorial}} as follows:
\begin{equation}\label{16}
(x)_{n,\lambda}=x(x-\lambda)(x-2\lambda)\cdots(x-(n-1)\lambda).
\end{equation}
Note that $(x)_{n,1}=x(x-1)\cdots(x-n+1)=(x)_n$, $n\geq0$. Therefore, by \eqref{9} and \eqref{16}, we obtain the following theorem.
\begin{theorem}\label{thm3}
For $n \geq 0$, we have
\begin{equation*}
\iz\left([x+y]_q\right)_{n,\lambda}d\mu_{-q}(y)=E_{n,q}(x|\lambda),
\end{equation*}
where $\left([x+y]_q\right)_{n,\lambda}=[x+y]_q\left([x+y]_q-\lambda\right)\cdots\left([x+y]_q-(n-1)\lambda\right)$.
\end{theorem}
From \eqref{7}, we can easily derive the following equation:
\begin{equation}\label{17}
\begin{split}
\iz[x+y]_q ^nd\mu_{-q}(y)=&\frac{[2]_q}{(1-q)^n}\sum_{l=0} ^n \binom{n}{l}(-1)^lq^{lx}\frac{1}{1+q^{l+1}}\\
=&[2]_q\sum_{m=0} ^{\infty}(-1)^m[x+m+1]_q ^n.
\end{split}
\end{equation}
Thus, by \eqref{17}, we get the generating function of the $q$-Euler polynomials as follows:
\begin{equation}\label{18}
\sum_{n=0} ^{\infty}E_{n,q}(x)\frac{t^n}{n!}=[2]_q\sum_{m=0} ^{\infty}(1-)^me^{[x+m+1]_q}.
\end{equation}
From \eqref{9} and \eqref{18}, we can derive the generating function of degenerate $q$-Euler polynomials which is given by
\begin{equation}\label{19}
\sum_{n=0} ^{\infty}E_{n,q}(x|\lambda)\frac{t^n}{n!}=[2]_q\sum_{m=0} ^{\infty}(-1)^m\left(1+\lambda t\right)^{\frac{[x+m+1]_q}{\lambda}}.
\end{equation}
By \eqref{19}, we easily get
\begin{equation}\label{20}
\sum_{n=0} ^{\infty}E_{n,q}(x|\lambda)\frac{t^n}{n!}=\sum_{n=0} ^{\infty}\left([2]_q\sum_{m=0} ^{\infty}(-1)^m\left([x+m+1]_q\right)_{n,\lambda}\right)\frac{t^n}{n!}.
\end{equation}
Therefore, by \eqref{20}, we obtain the following theorem.
\begin{theorem}\label{thm4}
For $n \geq 0$, we have
\begin{equation*}
E_{n,q}(x|\lambda)=[2]_q\sum_{m=0} ^{\infty}(-1)^m\left([x+m+1]_q\right)_{n,\lambda}.
\end{equation*}
\end{theorem}

We observe that
\begin{equation}\label{21}
\begin{split}
(1+\lambda t)^{\frac{[x+y]_q}{\lambda}}=&(1+\lambda t)^{\frac{[x]_q}{\lambda}}(1+\lambda t)^{\frac{q^x[y]_q}{\lambda}}\\
=&\left(\sum_{m=0} ^{\infty}\left([x]_q\right)_{m,\lambda}\frac{t^m}{m!}\right)\left(\sum_{l=0} ^{\infty}\lambda^{-l}q^{lx}\frac{[y]_q ^l \left(\log(1+\lambda t)\right)^l}{l!}\right)\\
=&\sum_{n=0} ^{\infty}\left(\sum_{k=0} ^n \sum_{l=0} ^k\left([x]_q\right)_{n-k,\lambda}\lambda^{k-l}q^{lx}[y]_q ^lS_1(k,l)\binom{n}{k}\right)\frac{t^n}{n!}.
\end{split}
\end{equation}
From \eqref{7} and \eqref{21}, we can derive the following theorem.
\begin{theorem}\label{thm5}
For $n \geq 0$, we have
\begin{equation*}
E_{n,q}(x|\lambda)=\sum_{k=0} ^n \sum_{l=0} ^k\binom{n}{k}\left([x]_q\right)_{n-k,\lambda}\lambda^{k-l}q^{lx}S_1(k,l)E_{l,q}.
\end{equation*}
\end{theorem}

It is not difficult to show that
\begin{equation}\label{22}
q\iz\left([x+y+1]_q\right)_{n,\lambda}d\mu_{-q}(y)+\iz\left([x+y]_q\right)_{n,\lambda}d\mu_{-q}(y)=[2]_q\left([x]_q\right)_{n,\lambda}.
\end{equation}
Thus, by Theorem \ref{thm5} and \eqref{22}, we obtain the following theorem.
\begin{theorem}\label{thm6}
For $n \geq 0$, we have
\begin{equation*}
qE_{n,q}(x+1|\lambda)+E_{n,q}(x|\lambda)=[2]_q\left([x]_q\right)_{n,\lambda}.
\end{equation*}
\end{theorem}

Let $r\in{\mathbb{N}}$. Now, we recall the Carlitz's $q$-Euler polynomials of order $r$ which are given by the generating function to be
\begin{equation}\label{23}
\iz\cdots\iz e^{[x_1+\cdots+x_r+x]_q}d\mu_{-q}(x_1)\cdots d\mu_{-q}(x_r)\\
=\sum_{n=0} ^{\infty}E_{n,q} ^{(r)}(x)\frac{t^n}{n!}.
\end{equation}
Thus, by \eqref{23}, we get
\begin{equation}\label{24}
\iz\cdots\iz[x_1+\cdots+x_r+x]_q ^nd\mu_{-q}(x_1)\cdots d\mu_{-q}(x_r)=E_{n,q} ^{(r)}(x),~(n\geq0),{\text{ (see \cite{05})}}.
\end{equation}

From \eqref{24}, we have
\begin{equation}\label{25}
\begin{split}
E_{n,q} ^{(r)}(x)=&\frac{[2]_q ^r}{(1-q)^n}\sum_{l=0} ^n \binom{n}{l}(-1)^lq^{lx}\left(\frac{1}{1+q^{l+1}}\right)^r\\
=&[2]_q ^r\sum_{m=0} ^{\infty}(-1)^m\binom{r+m-1}{m}q^m[m+x]_q ^n.
\end{split}
\end{equation}
Thus, by \eqref{25}, we get the generating function of Carlitz's $q$-Euler polynomials of order $r$ which given by
\begin{equation}\label{26}
\begin{split}
\sum_{n=0} ^{\infty}E_{n,q} ^{(r)}(x)\frac{t^n}{n!}=[2]_q ^r\sum_{m=0} ^{\infty}(-1)^m\binom{r+m-1}{m}q^me^{[m+x]_q t}.
\end{split}
\end{equation}
We consider {\it{degenerate $q$-Euler polynomials of order $r$}} which are given by the generating function to be
\begin{equation}\label{27}
\begin{split}
\sum_{n=0} ^{\infty}E_{n,q} ^{(r)}(x|\lambda)\frac{t^n}{n!}=[2]_q ^r\sum_{m=0} ^{\infty}(-1)^m\binom{r+m-1}{m}q^m(1+\lambda t)^{[m+x]_q t}.
\end{split}
\end{equation}
Thus, by \eqref{27}, we have
\begin{equation}\label{28}
\begin{split}
&\sum_{n=0} ^{\infty}E_{n,q} ^{(r)}(x|\lambda)\frac{t^n}{n!}\\
=&[2]_q ^r \sum_{n=0} ^{\infty}\left(\sum_{m=0} ^{\infty}(-1)^m\binom{r+m-1}{m}q^m\left([m+x]_q\right)_{n,\lambda}\right)\frac{t^n}{n!}.
\end{split}
\end{equation}
By comparing the coefficients on the both sides of \eqref{28}, we obtain the following theorem.
\begin{theorem}\label{thm7}
For $n\geq 0$, we have
\begin{equation*}
E_{n,q} ^{(r)}(x|\lambda)=[2]_q ^r\sum_{m=0} ^{\infty}(-1)^m\binom{r+m-1}{m}q^m\left([m+x]_q\right)_{n,\lambda}.
\end{equation*}
\end{theorem}

From \eqref{5} and Theorem \ref{thm7}, we note that
\begin{equation}\label{29}
\sum_{n=0} ^{\infty}E_{n,q} ^{(r)}(x|\lambda)\frac{t^n}{n!}=\iz\cdots\iz(1+\lambda t)^{\frac{[x_1+\cdots+x_r+x]_q}{\lambda}}d\mu_{-q}(x_1)\cdots d\mu_{-q}(x_r).
\end{equation}
Thus, by \eqref{29}, we get
\begin{equation}\label{30}
\begin{split}
E_{n,q} ^{(r)}(x|\lambda)
=&\sum_{m=0} ^{\infty}\iz\cdots\iz\left([x_1+\cdots+x_r+x]_q\right)_{n,\lambda} d\mu_{-q}(x_1)\cdots d\mu_{-q}(x_r)\\
=&\sum_{m=0} ^nS_1(n,m)\lambda^{n-m}\iz\cdots\iz[x_1+\cdots+x_r+x]_q ^md\mu_{-q}(x_1)\cdots d\mu_{-q}(x_r).
\end{split}
\end{equation}
Therefore, by \eqref{23} and \eqref{30}, we obtain the following theorem.
\begin{theorem}\label{thm8}
For $n \geq 0$, we have
\begin{equation*}
E_{n,q} ^{(r)}(x|\lambda)=\iz\cdots\iz\left([x_1+\cdots+x_r+x]_q\right)_{n,\lambda} d\mu_{-q}(x_1)\cdots d\mu_{-q}(x_r).
\end{equation*}
Furthermore,
\begin{equation*}
E_{n,q} ^{(r)}(x|\lambda)=\sum_{m=0} ^n S_1(n,m)\lambda^{n-m}E_{m,q} ^{(r)}(x).
\end{equation*}
\end{theorem}

Replacing $t$ by $\frac{1}{\lambda}\left(e^{\lambda t}-1\right)$ in \eqref{29}, we get
\begin{equation}\label{31}
\begin{split}
&\iz\cdots\iz e^{[x_1+\cdots+x_r+x]_q t}d\mu_{-q}(x_1)\cdots d\mu_{-q}(x_r)\\
=&\sum_{m=0} ^{\infty}E_{m,q} ^{(r)}(x|\lambda)\frac{1}{m!}\lambda^{-m}\left(e^{\lambda t}-1\right)^m\\
=&\sum_{m=0} ^{\infty}E_{m,q} ^{(r)}(x|\lambda)\lambda^{n-m}S_2(n,m)\frac{t^n}{n!}\\
=&\sum_{n=0} ^{\infty}\left(\sum_{m=0} ^n \lambda^{n-m}E_{m,q} ^{(r)}(x|\lambda)S_2(n,m)\right)\frac{t^n}{n!}.
\end{split}
\end{equation}
Therefore, by \eqref{23} and \eqref{31}, we obtain the following theorem.
\begin{theorem}
For $n\geq 0$, we have
\begin{equation*}
E_{n,q} ^{(r)}(x)=\sum_{m=0} ^nE_{m,q} ^{(r)}(x|\lambda)\lambda^{n-m}S_2(n,m).
\end{equation*}
\end{theorem}

\end{document}